\documentclass[journal]{IEEEtran}

\usepackage{xcolor,soul,framed} 

\colorlet{shadecolor}{yellow}
\usepackage[pdftex]{graphicx}
\graphicspath{{../pdf/}{../jpeg/}}
\DeclareGraphicsExtensions{.pdf,.jpeg,.png}

\usepackage[cmex10]{amsmath}
\usepackage{array}
\usepackage{mdwmath}
\usepackage{mdwtab}
\usepackage{eqparbox}
\usepackage{hyperref}
\hypersetup{
     colorlinks =true,
     linkcolor = blue,
     filecolor = blue,
     citecolor = magenta,      
     urlcolor = magenta,
     }
\usepackage{xurl}
\usepackage{cleveref,bm}
\usepackage{float}
\usepackage{amssymb}
\usepackage{algorithm}
\usepackage[noend]{algpseudocode}
\usepackage{graphics}
\usepackage{url}

\usepackage[caption=false, font=footnotesize]{subfig}
\allowdisplaybreaks

\begin{document}
\bstctlcite{IEEEexample:BSTcontrol}
    \title{Heuristics for Multi-Vehicle Routing Problem Considering Human-Robot Interactions}
  \author{Venkata Sirimuvva Chirala,
      Kaarthik Sundar, Saravanan Venkatachalam,
      Jonathon M. Smereka,
      Sam Kassoumeh 

  \thanks{Venkata Sirimuvva Chirala and Saravanan Venkatachalam are with the Department of Industrial and Systems Engineering, Wayne State University, Detroit, MI, 48202 USA (e-mail: \texttt{go1577@wayne.edu}; \texttt{fz1533@wayne.edu}).}%
  \thanks{Kaarthik Sundar is with Los Alamos National Labaratory, Los Alamos, New Mexico, 87545 USA (e-mail: \texttt{kaarthik@lanl.gov}).}
  \thanks{Jonathon M. Smereka and Sam Kassoumeh are with US Army CCDC, Ground Vehicles Systems Center, MI, 48397 USA (e-mail: \texttt{jonathon.m.smereka.civ@army.mil}; \texttt{sam.a.kassoumeh.civ@army.mil}).}}  

\maketitle

\begin{abstract}
Unmanned ground vehicles (UGVs) are being used extensively in civilian and military applications for applications such as underground mining, nuclear plant operations, planetary exploration, intelligence, surveillance and reconnaissance (ISR) missions and manned-unmanned teaming. We consider a multi-objective, multiple-vehicle routing problem in which teams of manned ground vehicles (MGVs) and UGVs are deployed respectively in a leader-follower framework to execute missions with differing requirements for MGVs and UGVs while considering human-robot interactions (HRI). HRI studies highlight the costs of managing a team of follower UGVs by a leader MGV. This paper aims to compute feasible paths, replenishments, team compositions and number of MGV-UGV teams deployed such that the requirements for MGVs and UGVs for the missions are met and the path, replenishment, HRI and team deployment costs are at minimum. The problem is first modeled as a a mixed-integer linear program (MILP) that can be solved to optimality by off-the-shelf commercial solvers for small-sized instances. For larger instances, a variable neighborhood search algorithm is offered to compute near optimal solutions and address the challenges that arise when solving the combinatorial multi-objective routing optimization problem. Finally, computational experiments that corroborate the effectiveness of the proposed algorithms are presented.
\end{abstract}


\IEEEpeerreviewmaketitle

\section{Introduction}

Industry 4.0 \cite{bai2020industry} concepts and technologies such as artificial intelligence (AI), machine learning, cloud computing are establishing their presence across innumerable fields in the recent years. Worldwide artificial intelligence (AI) software revenue is forecast to total \$62.5 billion in 2022, an increase of 21.3\% from 2021 \cite{Gartner}. Advancement in AI technology is causing a paradigm shift in the way autonomous vehicle technology is applied across various civilian and military applications. UGVs and unmanned aerial vehicles (UAVs), previously controlled manually or remotely, with the incorporation of AI are able to operate with increasing autonomy. The use of UGVs and UAVs to assist humans can be seen across various industries such as agricultural spraying and harvesting, nuclear plant operations, crowd control, environmental monitoring, and ISR missions such as explosive ordnance disposal (EOD), firefighting, search and rescue.

\begin{figure}
	\centering
	\includegraphics[scale=0.45]{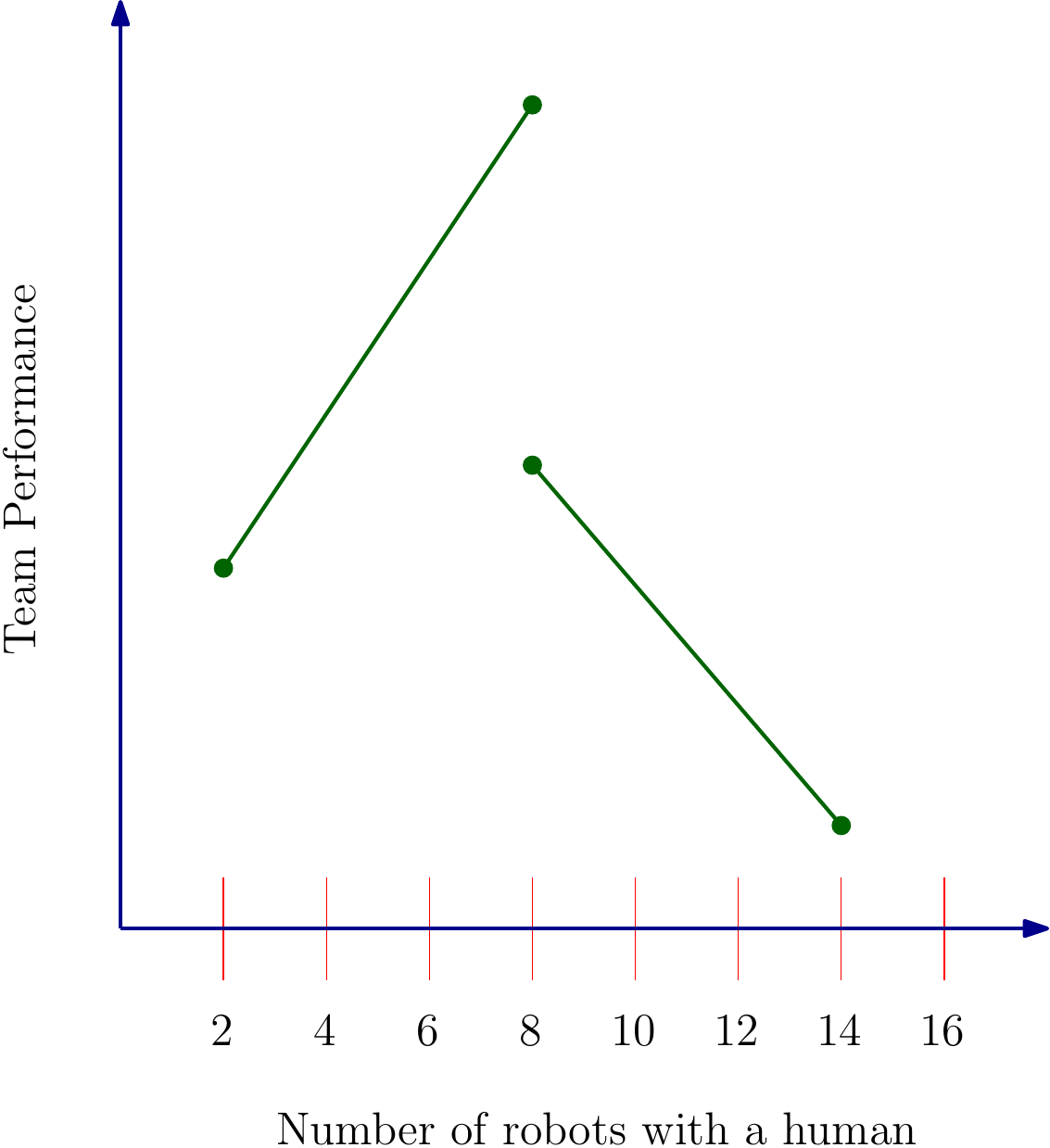}
	\caption{A piece-wise linear discontinuous function for the team's performance based on a HRI study in \cite{velagapudi2008scaling}}
	\label{FigHRI}
\end{figure}

In 2017, the United States Army published its strategy for the development and deployment of robotic and autonomous systems \cite{USArmyBook}. This strategy emphasizes primarily on human-machine collaboration by integrating army with robotic and autonomous systems solutions. Teaming humans/MGVs with robots (UGVS and UAVs) and AI will provide significant advantage to military forces by increasing combat effectiveness, lightening the soldiers' cognitive and physical workloads, facilitating maneuvering to difficult terrains and reducing harm to human personnel. MGV-UGV teams in a leader-follower framework are deployed in semi-autonomous convoy operations thereby providing logistics formations the ability to rest drivers for critical tasks only humans/MGVs can perform \cite{USArmyBook}. A convoy is a group of vehicles travelling together accompanied by armed troops or other vehicles. Within the leader-follower framework, the leader refers to an MGV. The MGV leader lays down ``virtual train tracks" which a set of follower UGVs use to follow the leader. Determining ideal team compositions, i.e., the ratio of human operators to robots is crucial for missions involving MGV-UGV teams \cite{yanco2002taxonomy}, \cite{liu2013robotic}. The team compositions used in this effort are based on HRI studies. HRI is a rapidly evolving field. The ideal team compositions for MGV-UGV teams are derived from various characteristics such as neglect tolerance and activity time. The autonomy of a robot is measured by its neglect time. A human operator can ignore a robot for the duration of its neglect time. Activity time measures the amount of time a human operator is involved in instructing a robot \cite{olsen2004fan}. These HRI studies provide a model for the human costs of managing a team of UGVs by MGVs \cite{goodrich2008human}. Authors in \cite{velagapudi2008scaling} conducted an urban search and rescue mission (USAR) to evaluate the task performance while a human operator was carrying varying number of robots. Fig. \ref{FigHRI} explains the HRI/team composition costs based on \cite{velagapudi2008scaling}. Authors in \cite{olsen2004fan} also measure task effectiveness as a function of the number of robots carried by a human operator. 

In this research, we consider a multiple vehicle routing problem (MVRP) in which MGV-UGV teams are deployed to implement missions requiring human-robot collaboration. Given a set of point of interests (POIs) with differing requirements for MGVs and UGVs, the objective is to find ideal number of teams, paths and replenishments for each team, and team compositions such that the following conditions are satisfied:

\begin{enumerate}
    \item The demand for the MGVs and UGVs at the POIs are satisfied. The demand for the UGVs is met either by the MGV-UGV team travelling to the POIs or by a combination of the former and additional replenishments supplied to the POI from the depot that then join the MGV-UGV team at the POI.
    \item The MGV-UGV teams follow a leader-follower framework where follower UGVs follow a leader MGV based on HRI studies with an exception for the replenishment UGVs supplied from the depot until they join the MGV-UGV team at a POI.
    \item The sum of travel (path and replenishment) costs, HRI (team composition) costs and costs for the number of teams deployed to complete the missions are minimum.
\end{enumerate}

\begin{figure}
    \centering
  \subfloat[\label{5a}]{%
       \includegraphics[scale = 0.55]{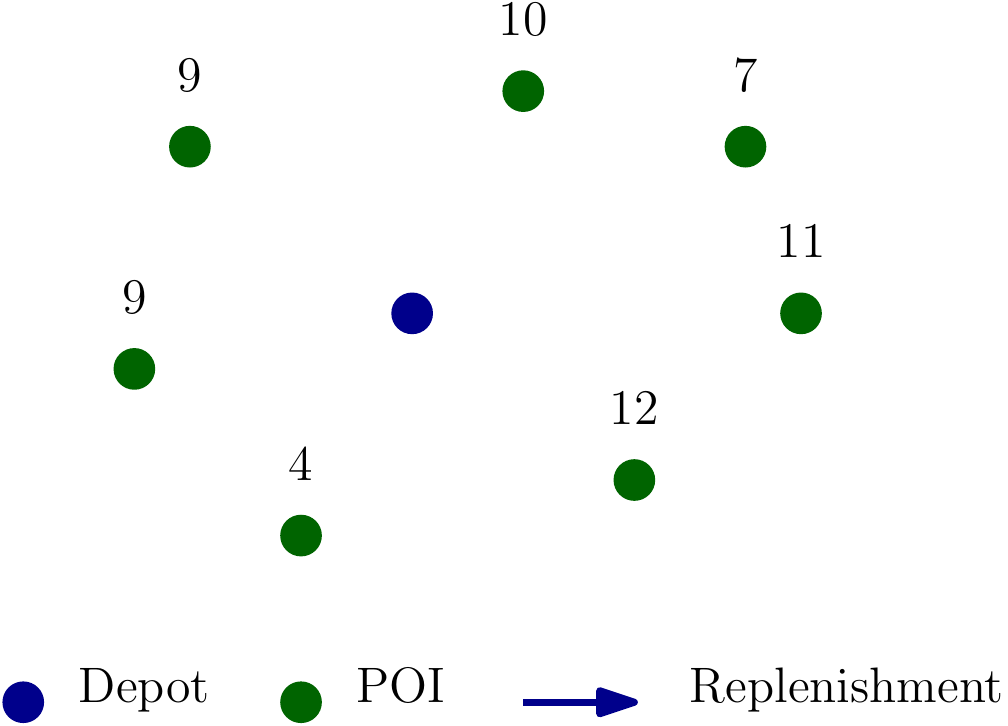}}
  \\
  \subfloat[\label{5b}]{%
        \includegraphics[scale = 0.55]{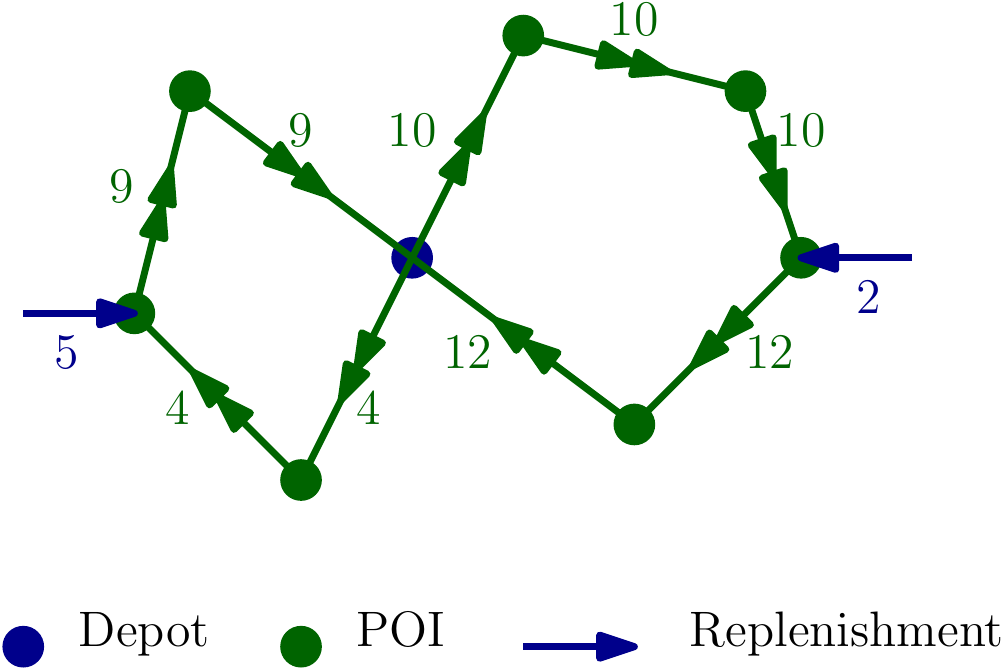}}
  \caption{ (a) shows the location of depot, POIs and the demand for UGVs at the POIs. (b) shows feasible routes such that all POIs are visited, and feasible assignments such that the demand for the UGVs at the POI is satisfied: replenishment UGVs from depot (blue arrows) and total UGVs traversing the edges (green arrows), given two MGV-UGV teams.}
  \label{feasible solution} 
\end{figure}

Given two vehicles, seven POIs and the demand for UGVs at POIs indicated in Fig \ref{5a}, a feasible solution for the MVRP instance is shown in Fig \ref{5b}. In the next section, we present a detailed overview of the literature that concerns this problem.


\section{Related Work}
The multiple vehicle, MGV-UGV team routing problem described in this paper is a combination of the capacitated vehicle routing problem (CVRP) \cite{ralphs2003capacitated} and the generalized assignment problem (GAP) \cite{cattrysse1992survey}. The CVRP problem determines the optimal routes for the MGV-UGV teams and the GAP problem is used to determine the optimal team compositions for each team; in this case, the costs for the team compositions are based on HRI studies.

The CVRP is a variant of the vehicle routing problem (VRP) \cite{toth2002vehicle} which determines the optimal routes for a fleet of capacitated vehicles initially located at one or more depots to serve a set of customers \cite{toth2002models}. There has been extensive research conducted in literature to apply exact methods \cite{LAPORTE1987147,doi:10.1057/palgrave.jors.2602345,fukasawa2006robust} that find the optimal solution to the CVRP using methods such as branch-and-cut. Heuristics \cite{szeto2011artificial,lin2009applying} that solve the problem faster and efficiently, however compromising on optimality, have also been explored. Given $n$ agents and $m$ tasks, the GAP determines optimal allocation of agents to tasks, such that each agent executes only one task subject to capacity restrictions on the agents \cite{cattrysse1992survey}. Similar to the CVRP, there is extensive research on exact methods \cite{nauss2003solving,park1998lagrangian,savelsbergh1997branch} and heuristics \cite{yagiura2004ejection,diaz2001tabu} for GAP.

In our earlier work presented in \cite{chirala2021multi}, we propose an optimization model to solve the MVRP problem. This was a first attempt at proposing a solution to multi-vehicle route planning problems involving MGV-UGV teams taking human-robot interaction into consideration. The optimization model is presented in this paper in Sec. \ref{sec:model} to keep the presentation self-contained. The CVRP problem and the GAP problem, the two key sub-problems of the MVRP, are NP-Hard \cite{lenstra1981complexity}. Therefore, scaling the optimization model to solve large scale instances is very challenging, time consuming and computationally expensive. Hence, the use of heuristic methods becomes inexorable.

VNS, proposed by Mladenovic and Hansen \cite{hansen2003tutorial}, is a very well-known metaheuristic used to solve such NP-hard optimization problems. VNS utilizes systematic change of neighborhoods within a local search algorithm exploring increasingly distant neighborhoods of the current incumbent solution and moving to a new solution only if an improvement is made \cite{mladenovic1997variable}. Authors in \cite{KYTOJOKI20072743} use the VNS with seven neighborhoods, three intra-route and four inter-route neighborhoods, alongside the guided local search metaheuristic to solve capacitated vehicle routing problems with a single depot and $2,400$-$20,000$ customers. They show that the proposed methodology is able to produce fast, high-quality feasible solutions for numerous large-scale benchmark instances within reasonable time. Authors in \cite{polacek2004variable} apply VNS algorithm to multi-depot vehicle routing problem with time windows for instances containing up to $300$ customers. They compare the performance of VNS to Tabu-search heuristic \cite{gendreau1994tabu}  and conclude that VNS is able to perform better with solution quality and scaling. VNS has also been an efficient heuristic for assignment problems \cite{mitrovic2009local, gavranovic2012variable}. Adaptation of VNS to multi-objective combinatorial optimization problems is gaining pace in the research community. Authors in \cite{duarte2015multi} suggest an adapted VNS approach to solve two combinatorial bi-objective optimization problems. They report that the suggested approach outperforms both Non-dominated Sorting Genetic Algorithm (NGSA) and the Strength Pareto Evolutionary Algorithm, two standard and popular algorithms used multi-objective optimization. Authors in \cite{HIDALGOPANIAGUA201620} apply the multi-objective VNS to a path planning problem in mobile robotics by optimizing three objectives, namely path safety to avoid obstacles, path length to control robot operation time and path smoothness to optimize energy consumption to obtain accurate and efficient paths. They apply VNS to eight different realistic scenarios and conclude that the multi objective VNS approach is superior in comparison to the NSGA. Given the efficacy of VNS in solving various classes of vehicle routing problems including problem with multiple vehicles, multiple depots, time windows, capacity restrictions and multi objective combinatorial optimization problems, in this research, we propose the variable neighborhood search (VNS) algorithm approach to the MVRP problem and address the challenge of scaling the approach to large test instances.

The remainder of the paper is organized as follows. The optimization model used to solve the MVRP problem is presented in section \ref{sec:model}. The variable neighborhood search algorithm and approach are presented in section \ref{sec:VNS}. Computational results for randomly generated instances are presented in section \ref{sec:results}. Finally, the paper concludes in section \ref{sec:conclusion}.


\section{Optimization Model} \label{sec:model}

The optimization model presented in our previous work in \cite{chirala2021multi} is included in this section with additional explanation for completeness.

\subsection{Problem Description} \label{subsec:problemdesc}

Let $P= \{p_1,\dots,p_l\}$ denotes the set of POIs. POIs are the locations at which MGV-UGV teams are required to execute missions requiring human-robot collaboration. Let $d_0$ denote the depot where $K$ different categories of heterogeneous MGVs and UGVs represented as $M_k$, $m \in M_k$, are initially stationed. $V = P\cup d_0$ denotes the set of vertices and $E$ denotes the set of edges joining any pair of vertices. For each edge $(i,j) \in E$, $c_{ij}$ represents the travel cost incurred while traversing the edge $(i,j)$. Set $N$ (indexed by $n$) is used to represent the number of UGVs in a team, and the parameter $h_n$ represents the HRI costs of managing `$n$' UGVs. Parameter $T$ represents the cost incurred for each MGV-UGV team deployed to execute the missions. For any set $S \subset V$, $\delta^+(S)=\{(i,j) \in E: i\in S, j\notin S\}$ and $\delta^-(S)=\{(i,j) \in E: i\notin S, j\in S\}$ denote the set of outgoing and incoming edges from the set $S$. When $S = \{i\}$, we shall simply write $\delta^+(i)$ and $\delta^-(i)$ instead of $\delta^+(\{i\})$ and $\delta^-(\{i\})$, respectively. Also, for any $A \subseteq E$, let $x(A) = \sum_{(i,j)\in A} x_{ij}$.

The decision variables for the MVRP are as follows: 
each edge $(i,j)\in E$ is associated with a variable $x_{ij}^m$ that equals $1$ if the edge $(i,j)$ is traversed by the `$m$\textsuperscript{th}' MGV, and $0$ otherwise. Variable $y_j^m$ represents the number of UGVs visiting the POI $j\in P$ along with the MGV $m$. Variable $\bar y_{nj}^m$ represents variable $y_j^m$ as a series of binary variables. Variable $z_j^m$ represents the number of UGV's provided as reinforcement to the POI $j\in P$ from the base depot $d_0$ for the MGV $m$. Variable $r^m_{nj}$ represents variable $z_j^m$ as a series of binary variables. Also, variable  $w_{ij}^m$ represents the number of UGVs that visit a POI $i$ then a POI $j$, where $i,j\in P$ along with the MGV $m$.

\subsection{Objective Function} \label{subsec:objfn}

The objective of the problem is to minimize the sum of (1) routing costs, (2) replenishment costs, (3) HRI costs and finally (4) the cost incurred by the teams to the mission. The expressions for each of these costs functions are as follows:
\begin{subequations}
    \begin{flalign}
    \text{Routing costs: } & \mathcal R^1 \triangleq \sum\limits_{\substack{(i,j) \in E, n\in N \\ m \in M_k, k\in K}} c_{ij}x_{ij}^m \label{eq:routing-costs} \\
    \text{Replenishement costs: } & \mathcal R^2 \triangleq \sum\limits_{\substack{(i,j) \in E, n\in N \\ m \in M_k, k\in K}} c_{d_0j}{r}^m_{nj} \label{eq:replenishment-costs} \\
    \text{HRI costs: } & \mathcal H \triangleq \sum\limits_{\substack{j \in V, n\in N \\ m \in M_k, k\in K}} h_n \bar{y}_{nj}^m \label{eq:hri-costs} \\ 
    \text{Team costs: } & \mathcal T \triangleq \sum\limits_{\substack{j \in V, \\ m \in M_k, k\in K}} T \cdot x_{0j}^m \label{eq:team-costs}
    \end{flalign}
    \label{eq:costs}
\end{subequations}

The objective function is the sum of the of above costs weighted by user provided parameters $\alpha$, $\beta$, and $\gamma$ as follows:

\begin{flalign}
\min ~ \alpha \cdot (\mathcal R^1 + \mathcal R^2) + \beta \cdot \mathcal H + \gamma \cdot \mathcal T \label{eq:obj2-1}
\end{flalign}

The weighted sum method is a classical approach used to solve multi-objective optimization problems. The optimal solution to the problem depends on the weights assigned to the objectives. Weights are assigned to the objective in proportion to the preference factor assigned to the objective \cite{deb2014multi}. Alternatively, a problem instance can be evaluated for different weights to give decision makers multiple perspectives.




\subsection{Routing Constraints} \label{subsec:constraints}

Routing constraints impose the requirements that all POIs $i\in P$ should be visited and that the demand for the number of MGVs and UGVs at the POIs is satisfied. Constraint \eqref{f1_1} ensures that a feasible solution is connected and are referred to as sub tour elimination constraints. The optimization problem is initially solved by relaxing these constraints. Once a feasible solution is obtained, the \textit{callback} feature offered by solvers such as CPLEX, Gurobi etc. are used to add the constraints to the model if the constraints in \eqref{f1_1} are violated. A constraint is said to be in violation if the tour does not contain the depot $d_0$. Constraint \eqref{f1_2} states that the number of MGV-UGV teams leaving the depot must come back to the depot after completing the missions at the visited POIs. Constraint \eqref{f1_3} states that and MGV-UGV team entering a POI must leave the POI. For each POI $i$, the pair of constraints in \eqref{f1_4} and \eqref{f1_5} state that the demand for MGVs and UGVs at POI $i$ are satisfied.

\begin{subequations}
		\begin{align}
		\begin{split}
		&x^m(S) \leqslant |S| - 1 \quad \forall (i,j) \in S, \\
		&S\subset V\setminus\{d_0\} : S \neq \emptyset, \quad \forall m \in M_k, k \in K, \end{split} \label{f1_1} & \\
		& x^m(\delta^+(d_0)) = x^m(\delta^-(d_0))  \quad \forall m \in M_k, k \in K, \label{f1_2} &\\
		& x^m(\delta^+(i)) = x^m(\delta^-(i))  \quad \forall  i \in P, m \in M_k, k \in K, \label{f1_3} &\\
		\begin{split}
		&\sum_{m \in M_k} x^m(\delta^+(i)) \geq a_{ik} \text{ and } \\
		&\sum_{m \in M_k} x^m(\delta^-(i)) \geq a_{ik} \quad \forall  i \in P, k \in K, \end{split}
		\label{f1_4} &\\
		&\sum_{m \in M_k} y^m_j \geq b_{jk} \quad \forall  j \in P, k \in K, \label{f1_5} &
		\end{align}
		\label{eq:routing}
	\end{subequations}
	
\subsection{HRI/Teaming Constraints} \label{subsec:teamingconstraints}

The teaming constraints decide the team sizes for the MGV-UGV teams. For each POI, constraint \eqref{f1_6} states that the number of UGV's along side MGV $m$ visiting the POI should be equal to the the incoming UGV's from the previously visited POI and the reinforcement provided to the current POI from the base depot. Constraints in \eqref{f1_7}-\eqref{f1_9} state that UGVs' travel alongside MGV $m$ if and only if the edge $x_{ij}^m$ is traversed. These constraints are used to convert the non-linear relationship ($\sum_{m \in M_k} w^m_{ij} = \sum_{m \in M_k} y^m_i$ only if $x_{ij}^m$ = 1) to a linear one. Constraint \eqref{f1_10} state the dependence of UGVs on an MGV following a leader-follower framework, and constraints \eqref{f1_11} define the capacity for UGVs. Constraints \eqref{f1_12} state that team size does not change for each manned GV `$m$', and the constraints \eqref{f1_13} are the special order sets of constraints (type 1) stating that only one binary variable can be chosen. Constraint \eqref{f1_14} is used to account for replenishment costs at each POI. Finally, constraint \eqref{f1_15} imposes restrictions on the decision variables.

\begin{subequations}
    \begin{align}
		&\sum_{m \in M_k, i \in P,d_0} w^m_{ij} + z^m_j = \sum_{m \in M_k} y^m_j \quad \forall j \in (P,d_0) \label{f1_6}&\\
		\begin{split}
		&\sum_{m \in M_k} w^m_{ij} \leq \sum_{m \in M_k} B*x^m_{ij} \\
		&\quad \forall i \in (P,d_0), \forall  j \in (P,d_0), k \in K, \end{split}
		\label{f1_7} &\\
		\begin{split}
		&\sum_{m \in M_k} w^m_{ij} \leq \sum_{m \in M_k} y^m_i \\ 
		&\quad \forall  j \in (P,d_0), \forall i \in (P,d_0), k \in K, \end{split} 
		\label{f1_8} &\\
		\begin{split}
		&\sum_{m \in M_k} w^m_{ij} \geq \sum_{m \in M_k} (y^m_i - (1 - x^m_{ij})*B) \\
		&\quad \forall i \in (P,d_0), \forall  j \in (P,d_0), k \in K, \end{split}
		\label{f1_9} &\\
		&0 \leqslant y^m_i \leqslant \ell_k x_{ji}^m \quad \forall (i,j) \in E, m \in M_k, k \in K,  \label{f1_10} & \\		
		&0 \leqslant \sum_{m\in M_k} y^m_i \leqslant \ell_k \quad \forall  k \in K,  \label{f1_11} & \\									
		&y^m_j = \sum_{n\in N} n_n\bar{y}_{nj}^m \quad \forall  j \in P, m \in M_k, k \in K,  \label{f1_12} & \\
		&\sum_{n\in N} \bar{y}^m_{nj} \leq 1 \quad \forall  m \in M_k, k \in K,  \label{f1_13} & \\
		&z^m_j = \sum_{n\in N} n_n{r}_{nj}^m \quad \forall  j \in P, m \in M_k, k \in K,  \label{f1_14} & \\			
		\begin{split}
		&x_{ij}^m \in \{0,1\}, y^m_i \in \mathbb{Z}^+, w^m_{ij} \in \mathbb{Z}^+, z^m_i \in \mathbb{Z}^+, \\
		&\bar{y}^m_{nj}, \bar{r}^m_{nj} \in \{0,1\}, 0 \leqslant s_{d}^m \leqslant 1, \\
		&\quad \forall  (i,j) \in E, n \in N, m \in M_k, k \in K. \end{split} \label{f1_15} &
		\end{align}
		\label{eq:teaming}
	\end{subequations}

The optimization model proposed to solve the MVRP contains the CVRP and GAP, both of which individually are combinatorial optimization problems. They are NP-Hard \cite{lenstra1981complexity}. Therefore, the optimization model also is combinatorial. The model can find optimal solution to small instances, the results of which are presented in \cite{chirala2021multi}. Scaling the optimization model to larger instances is computationally challenging. In the next section, we present the variable neighborhood search algorithm to address this challenge and compute high quality sub-optimal solutions for the MVRP.


\section{Variable Neighborhood Search} \label{sec:VNS}

The skewed variable neighborhood search algorithm (SVNS) is a variant of the VNS meta-heuristic that allows for the exploration of valleys far from the incumbent solution \cite{hansen2003tutorial}. SVNS is particularly advantageous when there may be near optimal solutions in separated and distant valleys. The exploration of far away valleys is done by allowing for \textit{`recentering'} of the search. A solution outside of the incumbent solution is accepted if it lies far away from the incumbent solution. \textit{Recentering} requires a metric to define the distance between the incumbent and current solutions. In this research, we use the SVNS, a variant of VNS.

\subsection{Construction Heuristic} \label{subsec:constheur}

SVNS starts with an initial feasible solution. Hence, we develop an initial feasible solution to the MVRP. Construction heuristic start with an empty solution to a given problem and extend the current solution repeatedly until a complete solution is found \cite{Wikipediaconst}. The initial feasible solution to the MVRP consists of two components: feasible routes i.e., the sequence of visits to the POIs by MGV-UGV teams such that all POIs are visited and feasible assignments i.e., the team compositions for MGV-UGV teams for each route. Feasible assignments ensure that the demand for the UGVs at the POI is met either by the incoming team to the POI from the previously visited POI or by the replenishment UGVs supplied from the base depot. We first construct the initial routes to the multiple travelling salesman problem (mTSP) by using an extension of the Lin-Kernighan-Helsgaun heuristic \cite{helsgaun2000effective}, LKH-3 heuristic solver \cite{LKH-3}. We then construct feasible assignments to the MVRP by comparing the following two assignment strategies and implementing the strategy that provides a better solution: 

\begin{enumerate}
   
    \item asgn-rule1: For each MGV-UGV team visiting a subset of POIs, consider the maximum UGV demand from the subset and provide this UGV demand from the depot allowing for no replenishments from the depot.
    \item asgn-rule2: For each MGV-UGV team visiting a subset of POIs, provide required UGV demand and replenish POIs as and when required based on the UGV demand at POIs.
    
\end{enumerate}

\begin{figure}
    \centering
  \subfloat[\label{1a}]{%
       \includegraphics[width=0.50\linewidth]{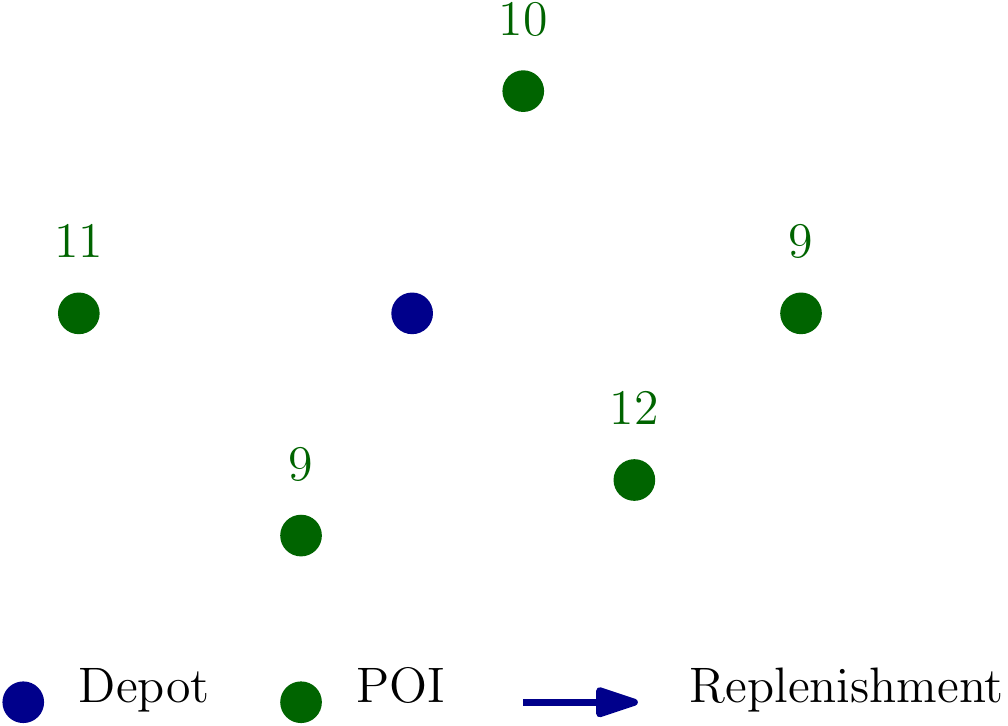}}
    \hfill
  \subfloat[\label{1b}]{%
        \includegraphics[width=0.48\linewidth]{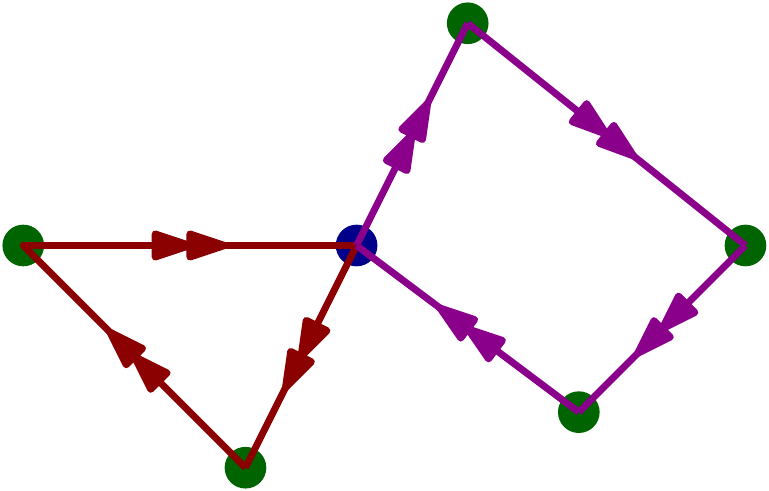}}
    \\
  \subfloat[\label{1c}]{%
        \includegraphics[width=0.48\linewidth]{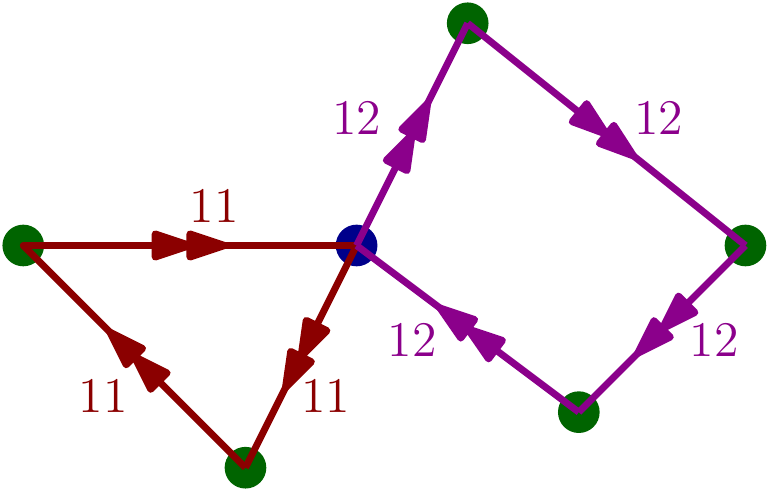}}
    \hfill
  \subfloat[\label{1d}]{%
        \includegraphics[width=0.48\linewidth]{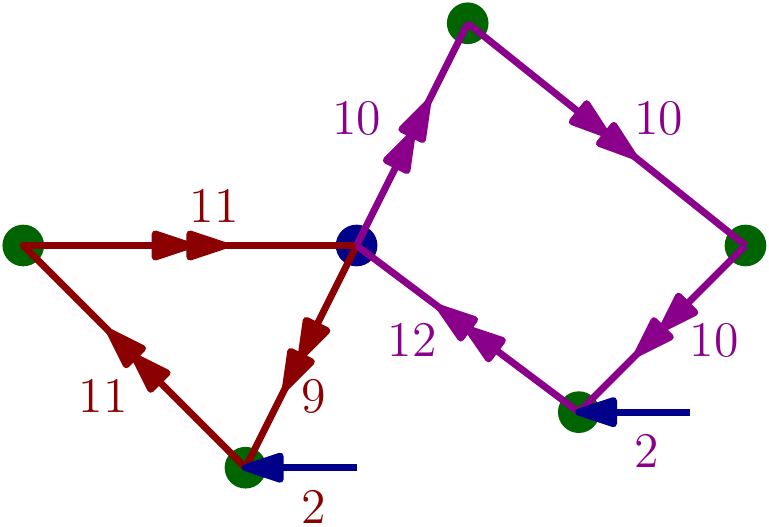}}
  \caption{ (a) shows the location of depot, POIs and the demand for UGVs at the POIs for MGV-UGV team missions, (b) shows a feasible routes with 2 MGV-UGV teams, (c) shows asgn-rule1 implementation, and (d) shows asgn-rule2 implementation.}
  \label{AssignmentStrategy} 
\end{figure}

The choice of assignment strategy implemented in the solution is dependent on the weights $\alpha$ and $\beta$ used in the objective function for the MVRP problem. asgn-rule1 strategy would provide a better solution when weight $\alpha$ is prioritized over $\beta$, since no replenishments are allowed leading to a replenishment cost of 0, however, increasing the HRI cost due to the increase in team sizes. Alternatively, asgn-rule2 increases the replenishment costs thereby reducing the team sizes and the HRI costs proving effective when weight $\beta$ is prioritized over $\alpha$. Figure \ref{AssignmentStrategy} explains the assignment strategies. Finally, we reverse the routes generated by the LHK-3 solver and construct feasible assignments once again as mentioned above. The best solution overall is used as the initial feasible solution for the SVNS. The costs for the construction heuristic solution may be high since the routes obtained from the LKH-3 heuristic may not be the best with respect to the assignments of the MGV-UGV teams. A different sequence of routes may yield a better suited solution for assignments. In the next section, we define the neighborhoods used in SVNS.

\subsection{Neighborhoods} \label{subsec:neighborhood}

\begin{figure}
    \centering
  \subfloat[\label{2a}]{%
       \includegraphics[width=0.45\linewidth]{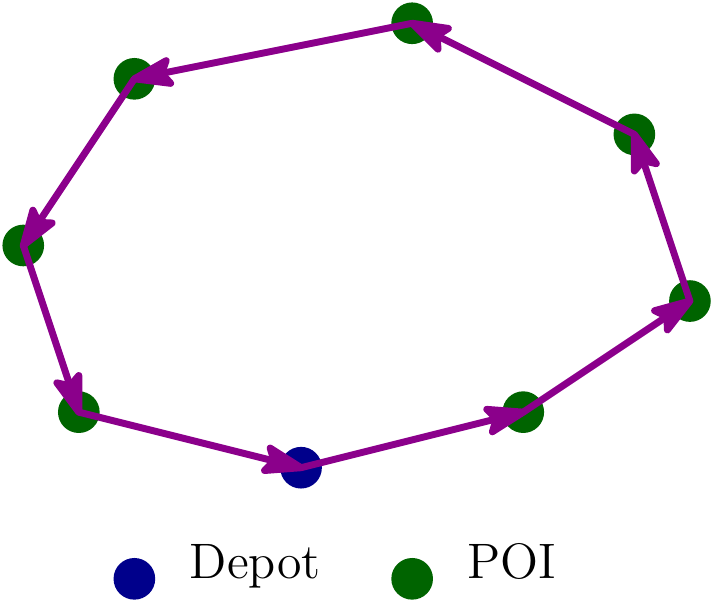}}
    \hfill
  \subfloat[\label{2b}]{%
        \includegraphics[width=0.45\linewidth]{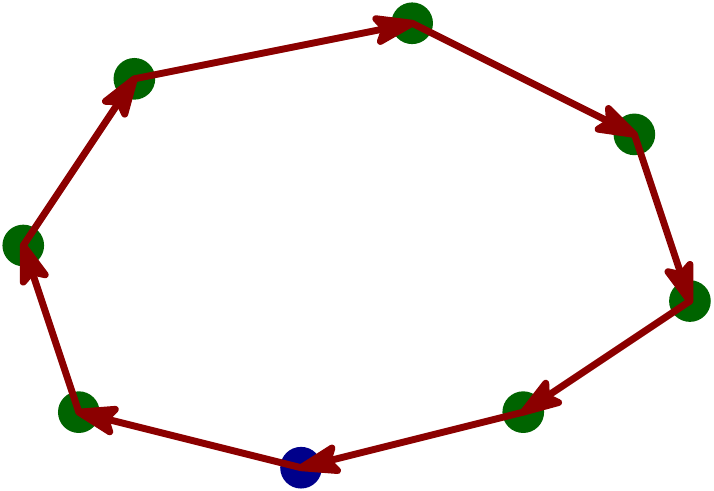}}
    \\
  \subfloat[\label{2c}]{%
        \includegraphics[width=0.45\linewidth]{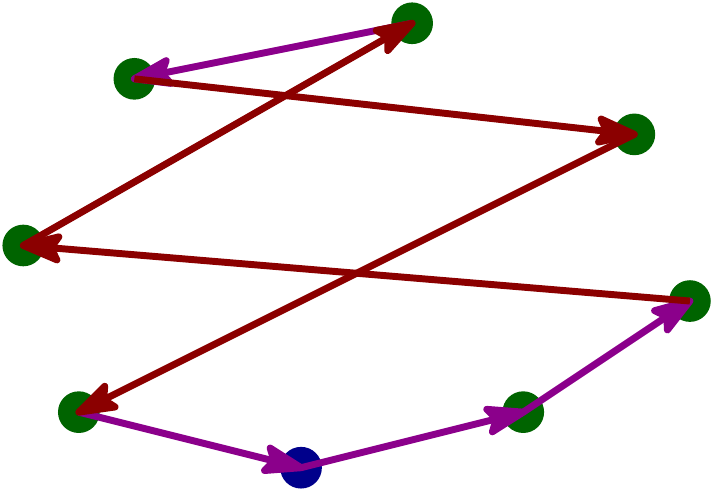}}
    \hfill
  \subfloat[\label{2d}]{%
        \includegraphics[width=0.45\linewidth]{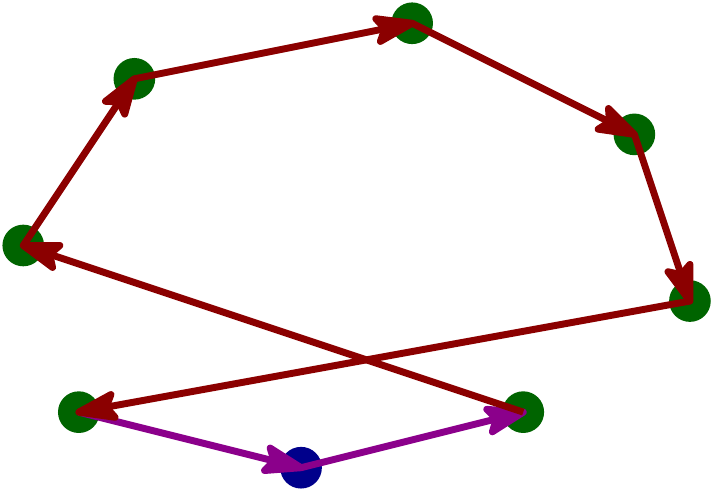}}
  \caption{ (a) shows the locations of the depot and POIs and path for a MGV-UGV team before applying intra neighborhood moves, (b) shows path after applying reverse neighborhood, (c) shows path after applying POI-swap-intra neighborhood, and (d) shows path after applying 2-opt-intra neighborhood.}
  \label{intra_route_neighborhoods} 
\end{figure}

\begin{figure}
    \centering
  \subfloat[\label{3a}]{%
       \includegraphics[width=0.45\linewidth]{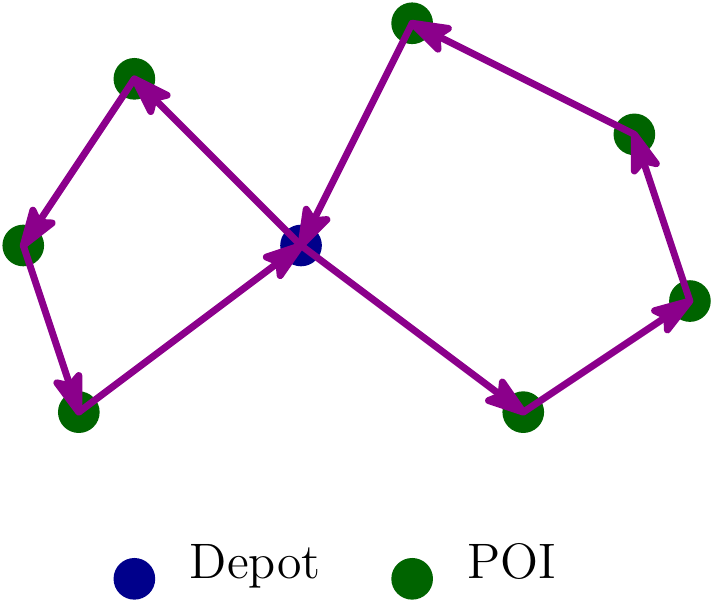}}
    \hfill
  \subfloat[\label{3b}]{%
        \includegraphics[width=0.45\linewidth]{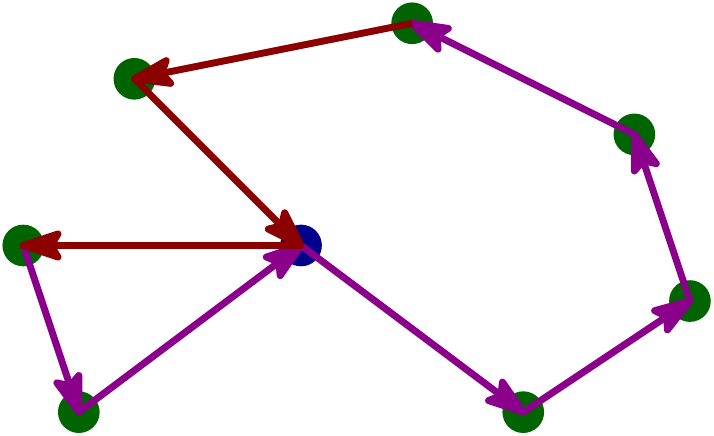}}
    \\
  \subfloat[\label{3c}]{%
        \includegraphics[width=0.45\linewidth]{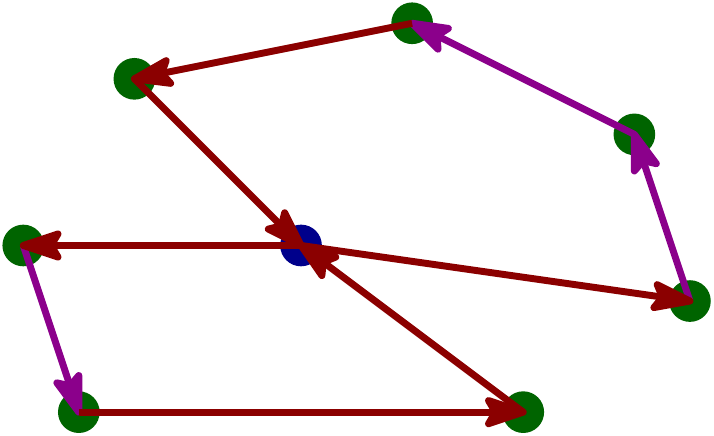}}
    \hfill
  \subfloat[\label{3d}]{%
        \includegraphics[width=0.45\linewidth]{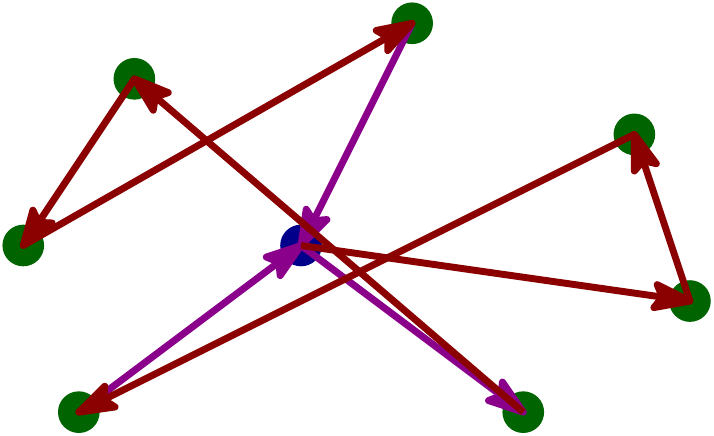}}
  \caption{ (a) shows the locations of the depot and POIs and paths for two MGV-UGV teams before applying inter neighborhood moves, (b) shows paths after applying POI-remove-insert neighborhood, (c) shows paths after applying POI-swap-inter neighborhood, and (d) shows paths after applying POI-SequenceExchange-inter neighborhood.}
  \label{inter_route_neighborhoods} 
\end{figure}

We define seven neighborhoods for the MVRP, three intra route neighborhoods, namely (i) reverse, (ii) POI-swap-intra, (iii) 2-opt-intra and three inter route neighborhoods, namely (iv) POI-remove-insert, (v) POI-swap-inter and (vi) POI-SequenceExchange-inter. Given routes for the MGV-UGV teams, the reverse neighborhood reverses the sequence of visits to POIs for each route. The POI-swap-intra neighborhood selects randomly, two POIs from a route and interchanges them with each other. 2-opt was first introduced by Croes\cite{croes1958method}. The 2-opt-intra neighborhood makes use of the 2-opt swap mechanism where, given a route and a segment of the route, the segment is reversed and reinserted back in its place. The intra neighborhoods are explained in Fig. \ref{intra_route_neighborhoods}. Given two routes, the POI-remove-insert takes a POI randomly from one route and inserts it into a random position in another route. The POI-swap-inter neighborhood takes a POI at random from two routes and exchanges the POIs. Finally, the POI-SequenceExchange-inter neighborhood takes one segment each, from two routes and and exchanges the segments. The inter route neighborhoods are explained in Fig. \ref{inter_route_neighborhoods}. For assignments, all intra route and inter route neighborhoods are evaluated with the two assignment strategies mentioned in section \ref{subsec:constheur} and the best solution is used. 

Additionally, we define an assignment neighborhood named (vii) asgn-rule3 as follows: Given a MGV-UGV team visiting a set of POIs, the asgn-rule3 neighborhood meets the UGV demand at the POIs by allowing replenishments only at a few selected POIs. For POIs where no replenishments are allowed the UGV demand is carried from the MGV-UGV team coming into the POI. This neighborhood allows for a balance of the replenishments and team sizes (replenishment and HRI costs) for the MVRP, in contrast to the two assignment strategy asgn-rule1 and asgn-rule2 in section \ref{subsec:constheur}.

\subsection{Algorithm} \label{subsec:algorithm}

The construction heuristic solution obtained in section \ref{subsec:constheur} is given as the initial feasible solution to SVNS. The SVNS algorithm consists of three steps: (i) shaking, (ii) local search, and (iii) recentering. The POI-swap-intra neighborhood presented in Sec. \ref{subsec:neighborhood} is used in the shaking phase of the algorithm. The shaking phase is used to escape from the incumbent solution. The 2-opt-intra neighborhood, all three inter route neighborhoods and the asgn-rule3 neighborhood are used in local search phase of the algorithm. The recentering step of the algorithm performs the following: if the solution from the shaking and local search phase $f^{''}$ is better than the incumbent solution computed thus far $f$, the search is recentered i.e., $f^{''}$ is assigned to $f$ and the new solution becomes the incumbent solution. This is done in line 15 of the SVNS algorithm. Additionally in the SVNS, line 23 of the algorithm, the search is also recentered when $f^{''}$ is worse than $f$, provided the distance between $f^{''}$ and $f$ is large enough \cite{9123228}. Given two feasible solutions where a solution represents the total cost for the MVRP, a solution $f^{''}$ is said to be better than $f$ if the total cost $f^{''}$ is strictly less than $f$ since the MVRP is a cost minimization problem. Solution $f^{''}$ is said to be worse than $f$ if the total cost $f^{''}$ is strictly greater than $f$. The distance metric for the MVRP is defined as the relative gap $(\%)$ in the total cost of the two solutions. Given two feasible solutions, $f^{''}$ and $f$, the distance between them is said to be large enough if the relative gap is greater than or equal to 20\% and less than or equal to 50\%. If the search is recentered, a parameter $k$ is initialized to 1, otherwise, the parameter $k$ is increased by 1 unit. Only when the total objective value decreases, the number of iterations with no improvement (UnImproved) is initialized to 0. In all the other cases, UnImproved is increased by 1. The pseudo-code of the SVNS is presented in Algorithm \ref{SVNS}.

\begin{algorithm}
\caption{Skewed Variable Neighborhood Search}
\label{SVNS}
\textbf{Input:} $N_j$ \textit{for Shaking phase:} POI-swap-intra \\
$N_j$ \textit{for Local Search Phase:} 2-opt-intra, POI-remove-insert, POI-swap-inter, POI-SequenceExchange-inter, asgn-rule3
\begin{algorithmic}[1]
    \State \textit{f} $\leftarrow$ construction heuristic solution.
    \While{UnImproved $<$ UnImproved\_max}
        \State $k$ $\leftarrow$ 0
        \While{$k$ $<$ $k\_max$}
            \State $f^{'}$ $\leftarrow$ Shake($f$) \Comment{Shake}
            \State $j$ $\leftarrow$ 0
            \While{j $<$ 5} \Comment{Local Search}
                \State Find best neighbor in $N_{j}(f^{'})$ $\rightarrow$ $f^{''}$
                \If{$f^{''}$ is better than $f^{'}$}
                    \State $f^{'}$ $=$ $f^{''}$
                    \State $j$ $\leftarrow$ 0
                \Else
                    \State $j$ $\leftarrow$ $j +$ 1
                    \State $f^{''}$ $=$ $f^{'}$
                \EndIf
                \If{$f^{''}$ is better than $f$} \Comment{Recenter}
                    \State $f$ $=$ $f^{''}$
                    \State $k$ $\leftarrow$ 0
                    \State UnImproved $\leftarrow$ 0
                    \If{$f$ is better than IncumbentSol}
                        \State IncumbentSol $\leftarrow$ $x$
                    \EndIf
                \Else
                    \State UnImproved$\leftarrow$UnImproved$+$1
                    \If{$f^{''}$ is worse than and far from $f$}
                        \State $f$ = $f^{''}$
                        \State $k$ $\leftarrow$ 0
                    \Else
                        \State $k$ $\leftarrow$ $k +$ 1
                    \EndIf
                \EndIf
            \EndWhile    
        \EndWhile   
    \EndWhile    
\end{algorithmic}
\textbf{Output: IncumbentSol}
\end{algorithm}


\section{Computational Results} \label{sec:results}

The SVNS algorithm for the MVRP is implemented using Python. All the computational experiments were run on a Intel(R) Core(TM) i5-6200U PC running at 2.30GHz and 4 GB RAM.

\subsection{Performance of SVNS algorithm}

To show the overall performance and effectiveness of the algorithm, we generate three sets of instances: 80 small scale instances with five POIs and two vehicles; 80 medium scale instances with 20 POIs, three and five vehicles; and 80 large scale instances with 40 POIs, six and eight vehicles. For each of these instances, the path traversal costs and the replenishment costs are measured as the Euclidean distance costs for the edges. The edge for replenishment cost is the edge connecting the depot to the POI where replenishment is supplied. The HRI (team compositions) are based on HRI studies. We choose HRI costs based on Fig \ref{FigHRI}. The costs for the number of teams deployed linearly increases as the number of vehicles deployed for the missions increase. The demand for the UGVs is uniformly distributed between one and 12 at the POIs. The total capacity for UGVs that can be carried by an MGV is 12. The weights for the objective function $\alpha$, $\beta$, and $\gamma$ can take values 0, 0.1, 0.3, 0.6, 1 such that the that the sum of weights adds up to one. 

\begin{figure}
	\centering
	\subfloat[\label{6a}]{%
	\includegraphics[scale=0.72]{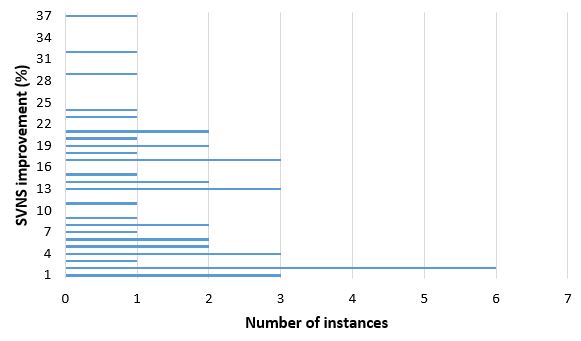}}
	\\
	\subfloat[\label{6b}]{%
	\includegraphics[scale=0.72]{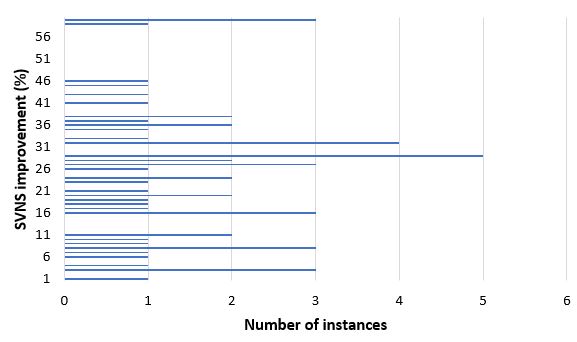}}
	\\
	\subfloat[\label{6c}]{%
	\includegraphics[scale=0.72]{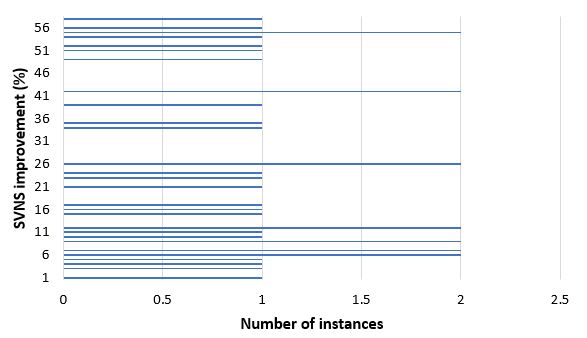}}
	\caption{Improvement of MVRP cost from construction heuristic (initial feasible solution) cost using SVNS for (a) small scale instances, (b) medium scale instances, (c) large scale instances.}
	\label{svns-imp}
\end{figure}

The small scale instance runs are made for both the optimization model presented in section \ref{sec:model} and the SVNS heuristic to evaluate the performance of SVNS. The average run time for small scale instances for the optimization model and SVNS are 4.28 and 2.02 seconds, respectively. SVNS results in optimal solutions for all 80 small scale instances. The initial feasible solution given by the construction heuristic in section \ref{subsec:constheur} alone produces optimal solutions for 38 instances of the 80 small scale instances. For the remaining 42 instances, the SVNS improves the initial feasible solution provided by the construction heuristic by a maximum of 36.35\% and an average improvement of 11.17\%. For the remainder of the paper, given 2 MVRP objective function costs, $f^{'}$ and $f^{''}$, improvement is defined as the decrease from $f^{'}$ to $f^{''}$ as a percentage of $f^{'}$. Fig \ref{6a} summarizes the performance of the SVNS algorithm for the remaining 42 small scale instances.

The optimization model due to being combinatorial in nature did not scale to find optimal solutions for medium and large scale instances. However with SVNS, medium scale instances take an average time of 11.6 seconds and 1.32 minutes for instances with three and five vehicles, respectively, and large scale instances take an average time of 5.79 and 13.21 minutes for instances with six and eight vehicles, respectively, to obtain sub-optimal solutions. SVNS improves the quality of the initial feasible solution by a maximum of 59.84\% and an average of 25.51\% for 56 of the 80 medium scale instances and a maximum of 57.34\% and an average of 25.48\% for 35 of the 80 large scale instances. Fig \ref{6b} and \ref{6c} summarize the performance of the SVNS algorithm for SVNS improved medium and large scale instances.

\begin{figure}
	\centering
	\includegraphics[scale=0.70]{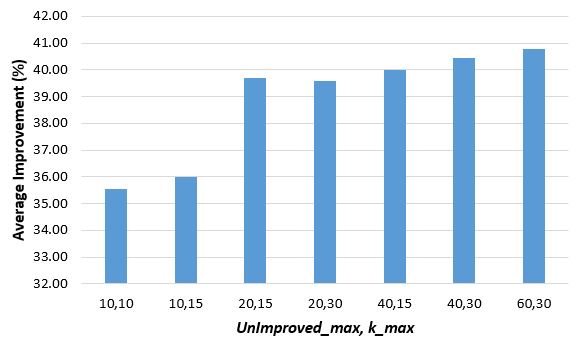}
	\caption{Improvement of MVRP cost from construction heuristic (initial feasible solution) cost using SVNS by varying parameters UnImproved\_max and $k\_max$ in Algorithm \ref{SVNS}}
	\label{parameter_fig}
\end{figure}

\subsection{Parametric Analysis}

Parameters, UnImproved\_max and $k\_max$ are used in line 2 and 4 respectively in the SVNS algorithm presented in Algorithm \ref{SVNS}. We generate 30 random instances with 10 small, 10 medium and 10 large scale instances to test the SVNS algorithm for different values of parameters, UnImproved\_max and $k\_max$. Figure \ref{parameter_fig} summarizes the average improvement in the MVRP cost from the construction heuristic using SVNS for differing values of the two parameters. The average improvement improves as the values of the parameters are increased. The lowest average improvement of 35.55\% occurs when both parameters UnImproved\_max and $k\_max$ are set to a value of 10. The average improvements are comparable when UnImproved\_max and $k\_max$ are set to 40 and 30 with 40.45\% (case 1) and 40 and 60 with 40.76\% (case 2) respectively. The time taken by SVNS increases by 4.01\% for the 10 large scale instances between case 1 and case 2. Based on the results, we set the parameter values of UnImproved\_max to 40 and $k\_max$ to 30. Increasing the values may yield better solutions, however, compromising on the time taken to solve the instances.

\begin{figure}
	\centering
	\includegraphics[scale=0.70]{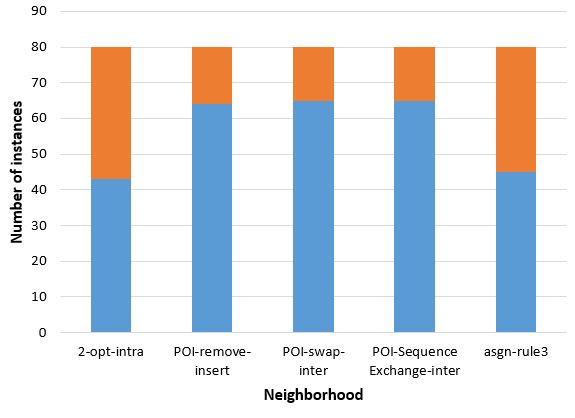}
	\caption{Instances resulting in optimal solution (indicated in blue color) from 80 small scale instances (indicated in orange color) when SVNS algorithm is used with one neighborhood in the local search phase at a time.}
	\label{neighborhood_optsol}
\end{figure}

\begin{table*}[t]
\centering
	\caption{SVNS improvement of MVRP cost from construction heuristic cost with individual neighborhoods in local search phase for small, medium and large scale instances }
	\label{tab:neighborhoods}
	\begin{tabular}{c c c c c c}
	& \\
		\hline
		\textit{Instance set} & \textit{Neighborhood} & \multicolumn{2}{c}{\textit{Number of instances improved}} & \textit{Maximum improvement} & \textit{Average improvement} \\
		 & &   \multicolumn{2}{c}{\textit{(Total instances: 80)}} & \textit{(\%)} & \textit{(\%)}    \\
		\hline
		
		Small scale& 2-opt-intra                & \multicolumn{2}{c}{13} & 28.56 & 10.79 \\
		& POI-remove-insert          & \multicolumn{2}{c}{37} & 31.81 & 11.73 \\
		& POI-swap-inter             & \multicolumn{2}{c}{37} & 36.35 & 11.27 \\
		& POI-SequenceExchange-inter & \multicolumn{2}{c}{37} & 36.35 & 11.27 \\
		& asgn-rule3                 & \multicolumn{2}{c}{14} & 28.56 & 9.77 \\
		
		Medium scale& 2-opt-intra                & \multicolumn{2}{c}{45} & 51.51 & 21.53 \\
		& POI-remove-insert          & \multicolumn{2}{c}{48} & 59.84 & 26.22 \\
		& POI-swap-inter             & \multicolumn{2}{c}{47} & 59.84 & 27.28 \\
		& POI-SequenceExchange-inter & \multicolumn{2}{c}{50} & 61.35 & 24.39 \\
		& asgn-rule3                 & \multicolumn{2}{c}{46} & 49.22 & 19.53 \\
		
		Large scale& 2-opt-intra                & \multicolumn{2}{c}{31} & 43.74 & 18.77 \\
		& POI-remove-insert          & \multicolumn{2}{c}{30} & 54.15 & 25.14 \\
		& POI-swap-inter             & \multicolumn{2}{c}{31} & 55.66 & 24.44 \\
		& POI-SequenceExchange-inter & \multicolumn{2}{c}{30} & 54.89 & 24.78 \\
		& asgn-rule3                 & \multicolumn{2}{c}{32} & 41.65 & 18.21 \\

		\hline
	\end{tabular}
\end{table*}

\subsection{Neighborhood Study}

Figure \ref{neighborhood_optsol} summarizes the number of instances that result in optimal solution for small scale instances when SVNS algorithm is used with one neighborhood in the local search phase at a time. SVNS with POI-swap-inter and POI-SequenceExchange-inter neighborhoods result in optimal solutions for 65 of the 80 small scale instances. The POI-remove-insert results in optimal solutions for 64 instances. All 3 neighborhoods are inter route neighborhoods and perform comparatively. They are also better in comparison to the 2-opt-intra, intra route neighborhood.

The performance of the SVNS algorithm for small, medium and large scale instances for each neighborhood operator used in the local search phase of the algorithm is summarized in Table \ref{tab:neighborhoods}. The POI-remove-insert neighborhood gives the highest average improvement for small and large scale instances. For medium scale instances, the highest average improvement is observed with the POI-swap-inter neighborhood. The performance of all three inter route neighborhoods: POI-remove-insert, and POI-SequenceExchange-Inter used in the local search of the SVNS are nearly identical. All the inter route neighborhoods perform better than the 2-opt-intra route neighborhood with small, medium and large scale instances. Overall, results suggest that the inter route neighborhoods in the local search phase of the SVNS show superior performance when compared to intra route neighborhood. The performance of the asgn-rule3 neighborhood improves as the instance size increases. However, the average improvement seen in asgn-rule3 neighborhood remains low as the neighborhood does not perform any route modifications.

\subsection{Sensitivity of weights on multi-objective function}

\begin{figure}
	\centering
	\includegraphics[scale=0.70]{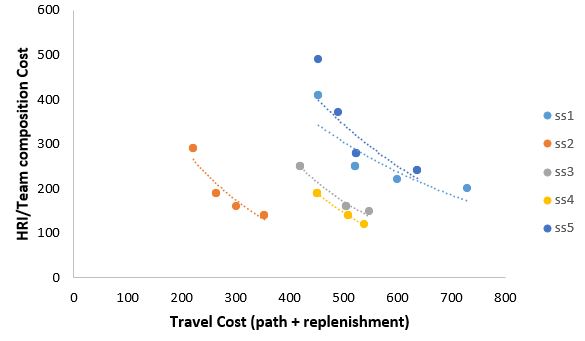}
	\caption{MVRP travel vs HRI cost plots for differing weights $\alpha$ and $\beta$ for five small scale instances (ss1, ss2, ss3, ss4 and ss5) to observe the effect of weights on the overall MVRP solution.}
	\label{Trade-off}
\end{figure}

\begin{table*}[t]
\centering
	\caption{Solutions for problem instances with differing weights ($\alpha$, $\beta$, $\gamma$)}
	\label{tab:solutionswithvaryingweights}
	\begin{tabular}{c c c c c c c c c c }
	& \\
		\hline
		\textit{Instance} & \multicolumn{3}{c}{\textit{User Parameters}} & \multicolumn{2}{c}{\textit{Travel Costs}} & \textit{HRI Cost} & \textit{Total Cost} & \multicolumn{2}{c}{\textit{Run Time (sec)}} \\
		& \textit{alpha}   &   \textit{beta}  &  \textit{gamma}  &  \textit{Path Cost}  &  \textit{Replenishment Cost} & & & \textit{SVNS} & \textit{opt}    \\
		\hline
		
		      ss1 & 1   &   1   &  1   & 521.89 & 0 & 250.00 & 871.89 & 1.43 & 6 \\
		      & 0.6 &   0.3 &  0.1 & 521.89 & 0 & 250.00 & 871.89 & 1.43 & 6    \\
		      & 0.6 &   0.1 &  0.3 & 454.15 & 0    & 410.00    & 964.15 & 1.84 & 5    \\
		      & 0.3 &   0.6 &  0.1 & 521.89 & 0 & 250.00 & 871.89 & 1.29 & 6    \\
		      & 0.1 &   0.6 &  0.3 & 533.59 & 196.34 & 200.03 & 1029.96 & 1.39 & 5    \\
		      & 0.1 &   0.3 &  0.6 & 510.47 & 89.56 & 220.03 & 920.06 & 1.26 & 6    \\
		      & 0.3 &   0.1 &  0.6 & 454.15 & 0    & 410.00    & 964.15 & 1.65 & 5    \\
		      
		      ss2 & 1   &   1   &  1   & 226.60 & 36.50 & 190.02 & 553.14 & 0.95 & 7 \\
		      & 0.6 &   0.3 &  0.1 & 226.60 & 36.50 & 190.02 & 553.14 & 0.63 & 7    \\
		      & 0.6 &   0.1 &  0.3 & 220.80 & 0    & 290.00    & 610.82 & 0.54 & 6    \\
		      & 0.3 &   0.6 &  0.1 & 264.80 & 36.50 & 160.03 & 561.36 & 0.54 & 5    \\
		      & 0.1 &   0.6 &  0.3 & 275.30 & 77.20 & 140.03 & 592.59 & 0.61 & 4    \\
		      & 0.1 &   0.3 &  0.6 & 275.30 & 77.20 & 140.03 & 592.59 & 0.49 & 6    \\
		      & 0.3 &   0.1 &  0.6 & 220.80 & 0    & 290.00    & 610.82 & 0.54 & 6    \\
		      
		      ss3 & 1   &   1   &  1   & 422.71 & 83.34 & 160.03 & 766.08 & 1.47 & 6 \\
		      & 0.6 &   0.3 &  0.1 & 420.82 & 0 & 250.00 & 770.82 & 1.27 & 5    \\
		      & 0.6 &   0.1 &  0.3 & 420.82 & 0    & 250.00    & 770.82 & 2.47 & 5    \\
		      & 0.3 &   0.6 &  0.1 & 422.71 & 83.34 & 160.03 & 766.08 & 1.31 & 4    \\
		      & 0.1 &   0.6 &  0.3 & 464.32 & 83.34 & 150.04 & 797.70 & 2.22 & 3    \\
		      & 0.1 &   0.3 &  0.6 & 422.71 & 83.34 & 160.03 & 766.08 & 1.24 & 6    \\
		      & 0.3 &   0.1 &  0.6 & 420.82 & 0    & 250.00    & 770.82 & 1.28 & 5    \\

		\hline
	\end{tabular}
\end{table*}

The weights for the objective function to the MVRP, given by $\alpha$, $\beta$, and $\gamma$ can either be predetermined by the decision makers to find the best solution based on the priority given to each of the objectives or pareto solutions obtained by differing the weights can be offered to the decision makers. In multi objective optimization (MOO), a solution is said to be non dominated, pareto optimal if none of the objective functions can be improved in value without degrading some of the other objective values \cite{MOO}. The weighted sum method is a classical \textit{A priori} method used to solve MOO problems. The set of objectives are scalarized to give a single objective by adding weights to each of the objectives. The weights $\alpha$, $\beta$, and $\gamma$ are parameters for scalarization. 

Fig \ref{Trade-off} plots the travel (path + replenishment) costs vs the HRI costs for differing weights $\alpha$, $\beta$, and $\gamma$ for a small-scale problem instance (ss1, ss2, ss3, ss4 and ss5) with five POIs and two vehicles. Table \ref{tab:solutionswithvaryingweights} details the complete solutions ss1, ss2 and ss3. A cost of 200 for each instance is added to the Total Cost in Table \ref{tab:solutionswithvaryingweights} to account for the team cost for deploying two vehicles. It can be observed that the travel costs (path and replenishments costs) are lesser when weight $\alpha$ is prioritized over $\beta$. This can be particularly advantageous when MGV-UGV teams are maneuvering through unfamiliar and difficult terrains and providing replenishments in such a scenario would be unsuitable and expensive. However, this increases the team sizes for the MGV-UGV teams, thereby increasing HRI/team composition costs. Alternatively, weight $\beta$ can be prioritized to reduce load on the MGVs and reduce HRI costs. However, such scenarios will have higher travel costs as observed. The best course of action for a given problem instance can be implemented based on the interests of the decision makers by studying the solution characteristics for the instance.


\section{Conclusion} \label{sec:conclusion}
In this research, we propose the SVNS algorithm to the multi-objective routing problem with multiple vehicles considering human robot interactions. Specifically, MGV-UGV teams are deployed in a leader-follower (MGV leader - UGV followers) framework to execute missions at POIs with differing requirements for MGVs and UGVs. To the best of our knowledge, our previous work in \cite{chirala2021multi} was a first attempt to model this problem. The model developed is a combinatorial optimization problem and is NP-Hard. SVNS is a variant of the VNS meta-heuristic that enhances the exploration of far away valleys. The SVNS addresses the challenges faced in scaling the previously modelled combinatorial optimization problem for the MVRP in \cite{chirala2021multi}. We generate three sets of instances: small, medium and large scale to evaluate the performance of SVNS. The SVNS produces optimal solutions for all small scale instances. For medium and large scale instances, SVNS generates high quality sub-optimal solutions for the optimization problem for more realistic and practical scenarios in reasonable time. Our future work involves focusing on extending the model to include uncertainties in the availability of UGVs, travel time, etc. Additional functionalities such as vehicle refueling can be incorporated into the model. Machine learning algorithms can be used to find the appropriate values for the weights, $\alpha$, $\beta$, and $\gamma$ used in the objective function.

\ifCLASSOPTIONcaptionsoff
  \newpage
\fi

\bibliographystyle{IEEEtran}
\bibliography{VNS}

\vfill

\end{document}